\documentclass[runningheads]{llncs}
\usepackage[T1]{fontenc}
\usepackage{graphicx}
\usepackage{amsmath, amssymb}  
\usepackage{booktabs}
\usepackage{comment}

\usepackage[ruled,vlined]{algorithm2e}
\usepackage{algpseudocode}

\SetAlgoLined
\DontPrintSemicolon
\SetKwInput{KwIn}{Input}
\SetKwInput{KwOut}{Output}
\SetKw{Return}{return}
\SetKwComment{Comment}{// }{}

\DeclareMathOperator{\diag}{diag}
\DeclareMathOperator{\conv}{conv}
\DeclareMathOperator{\rk}{rk}

\begin{document}
\title{Speeding up the Goemans-Williamson randomized procedure by difference-of-convex optimization}
\titlerunning{Speeding up the Goemans-Williamson procedure}
%
\author{Hadi Salloum\inst{1,2,3} \and
Roland Hildebrand\inst{1,2} \and
Nhat Trung Nguyen\inst{1} \and
Vitali Pirau\inst{1} \and
Amer Al Badr\inst{2,3} \and
Mohammad Alkousa\inst{2,3} \and
Alexander Gasnikov\inst{1,2}
}
\authorrunning{H. Salloum et al.}
%
\institute{MIPT, Institutskiy per. 9, 141701 Dolgoprudny, Russia \and
Innopolis University, Universitetskaya Str. 1, 420500 Innopolis, Russia
\email{\{salloum.h,khildebrand.r\}@mipt.ru; pireyvitalik@phystech.edu}\\
\and
Q Deep, Universitetskaya Str. 1, 420500 Innopolis, Russia\\
}

\maketitle              
\begin{abstract}
We present a novel approach to accelerate the Goemans–Wil\-liamson (GW) randomized rounding procedure for quadratic unconstrained binary optimization (QUBO) problems. Instead of solving the conventional semi-definite programming (SDP) relaxation, which is computationally expensive, we employ a difference-of-convex (DC) optimization framework to efficiently approximate the SDP solution. The DC optimization produces candidate vectors that are then used within the GW randomized rounding scheme to generate high-quality binary solutions. Furthermore, we perform direct expectation minimization over manifolds of matrices with limited rank to further enhance the solution quality. Our method is benchmarked on real-world QUBO instances, including inverse kinematics problems, and compared against state-of-the-art solvers, such as quantum-inspired algorithms, demonstrating competitive approximation guarantees alongside substantial computational gains.
\end{abstract}
\section{Introduction}

In this work we propose a new method for solving quadratic unconstrained binary optimization (QUBO) problems which minimize a quadratic function over the vertices of the hypercube. Our method, on the one hand, uses the classical ideas of Goemans, Williamson \cite{GoemansWilliamson}, and Nesterov \cite{Nesterov98} on the semi-definite relaxation of the QUBO problem and the associated randomized procedure for restoring a solution of the original problem, and, on the other hand, is inspired by the algorithm of Burer \cite{Burer02}, which searches over the manifold of factors of positive semi-definite matrices with rank not exceeding 2. The main idea of our algorithm is that it avoids solving the costly semi-definite program associated by the relaxation and merely searches over the manifold of factors of positive semi-definite matrices with rank not exceeding $k$, but in anticipation of the restoration of a solution to the original problem which is applied in the next step, minimizes not the extension of the original cost functional to this manifold, but the expectation of the value of the cost functional when applying the randomized procedure.

A strong motivation to formulate different discrete optimization problems as QUBO arose with the insight that looking for the minimum in a QUBO problem is equivalent to looking for the ground (minimum energy) state in an Ising model of coupled spins. Thus the solution of a QUBO can be realized by preparing a corresponding Ising state and letting it attain its ground state. Substantial effort has been spent to bring various problems in discrete optimization to QUBO form.

The QUBO formulation is somewhat restricted in its expressive power by the requirement that constraints on the involved binary variables have to be absent. This limits the scope of combinatorial problems that can be directly expressed as QUBO. OIn Section \ref{sec:applications} we present two such problems and the corresponding QUBO formulation. Note that these problems use $\{-1,+1\}$ binary variables. There are also QUBO formulations with $\{0,1\}$ variables, but these to can be converted into each other easily. 

Concrete problem classes which can be formalized as QUBO directly or after softening constraints are listed in Section \ref{sec:applications}. To this section we also delegate a description of the conversion between $\{0,1\}$ and $\{-1,+1\}$ QUBO problems. Other combinatorial problems can be brought to a QUBO form at the cost of softening constraints on the binary variables present in the original formulation. Among these are the quadratic assignment problem, vertex cover, set packing, and graph coloring.

An overview over existing methods to solve QUBO problems is delegated to the separate Section \ref{sec:QUBOmethods}. In an extra Section \ref{RelaxationBasedMethods} we present the semi-definite relaxation and the randomized procedure which are at the basis of our development.

\section{Applications and forms of QUBO} \label{sec:applications}

Since some of the problems are more easily reduced to a minimization problem over the vertices of the hypercube $[0,1]^n$, and others to a minimization problem over the vertices of the hypercube $[-1,+1]^n$, we shall first provide conversion formulas from one form of QUBO to the other. 

\medskip

Consider a QUBO of the form
\[ \min_{x \in \{0,1\}^n} (x^TAx + b^Tx) = \min_{x \in \{0,1\}^n} x^T(A + \diag(b))x,
\]
where the equality comes from the relation $x_i^2 = x_i$ for all involved $\{0,1\}$ binary variables. Let us substitute $y_i = 2x_i - 1 \in \{-1,+1\}$ and form a corresponding $\{-1,+1\}$ binary vector $y$. Then we have
\begin{align*}
    x^TAx + b^Tx &= \frac14(y + {\bf 1})^T(A + \diag(b))(y + {\bf 1}) \\ &= 
\frac14 \begin{pmatrix} 1 \\ y \end{pmatrix}^T \begin{pmatrix}
{\bf 1}^T(A + \diag(b)){\bf 1} & {\bf 1}^T(A + \diag(b)) \\
 (A + \diag(b)){\bf 1} & A + \diag(b) 
\end{pmatrix} \begin{pmatrix} 1 \\ y \end{pmatrix}.
\end{align*} 
Here ${\bf 1}$ stands for the all-ones vector. Therefore the original QUBO problem with $\{0,1\}$ binary variables is equivalent to the homogeneous problem
\[ \min_{z \in \{-1,+1\}^{n+1}}\,\frac14 \, z^T \begin{pmatrix}
                                                 {\bf 1}^T(A + \diag(b)){\bf 1} & {\bf 1}^T(A + \diag(b)) \\
                                                 (A + \diag(b)){\bf 1} & A + \diag(b) 
                                               \end{pmatrix} z
\]
with $\{-1,+1\}$ binary variables.

On the other hand, a homogeneous QUBO of the form
\[ \min_{y \in \{-1,+1\}^n}, y^TAy
\]
can be converted by the above substitution into the form
\begin{align*}
\min_{x \in \{0,1\}^n}\, (2x-{\bf 1})^TA(2x-{\bf 1}) &= \min_{x \in \{0,1\}^n}\,(4x^TAx-4{\bf 1}^TAx+{\bf 1}^TA{\bf 1}) \\ &= \min_{x \in \{0,1\}^n}\,4x^T(A - \diag(A{\bf 1}))x + {\bf 1}^TA{\bf 1}    
\end{align*} 
which is equivalent to a QUBO with $\{0,1\}$ binary variables.

We now have a choice to which form of QUBO to convert the considered combinatorial problem. We now consider two problems which can be equivalently reformulated as QUBO.

\medskip

\emph{Subset Sum:} This combinatorial problem can be formulated as follows. We are given positive integer numbers $w_1, \dots, w_n$ and have to find a subset $I \subset \{1, \dots, n\}$ such that
\[ \sum_{i \in I} w_i = \sum_{i \not\in I} w_i 
\]
or to show that no such subset exists. In other words, we look for a partition of the index set $\{1,\dots,n\}$ into subsets $I,\bar I$ such that the sum of the weights $w_i$ in each subset equals half the total weight. This problem is very hard, with the best-known algorithm proposing essentially an exhaustive search and having a time complexity of $O(2^{n/2})$ \cite{HorowitzSahni74}.

The problem can be formulated as a QUBO of the form
\[ \min_x x^TWx: \quad x \in \{-1,+1\}^n, 
\]
where $W = ww^T$ is a rank 1 matrix with factor $w = (w_1,\dots,w_n)^T \in \mathbb N_+^n$. 

If the optimal value of the problem is zero, then the solution $x^*$ yields the sought partition. If the optimal value is strictly positive, then no such partition exists.

\medskip

\emph{MaxCut:} Let $G = (V,E)$ be a given weighted undirected graph with edge weights $w_e > 0$ and vertex set $V = \{1,\dots,n\}$. Construct a symmetric matrix $W = (W_{ij})_{i,j=1,\dots,n}$ such that $W_{ij} = 0$ if $(i,j) \not\in E$ and $W_{ij} = w_e$ if $e = (i,j) \in E$. The goal is to partition the vertex set $V$ into subsets $I$ and $\bar I$ such that the cut weight 
\[ \sum_{i \in I, j \in \bar I} W_{ij} \]
is maximized.

Let us convert the MaxCut problem into a QUBO problem with $\{-1,+1\}$ binary variables. Represent the partition of $V$ by a $\{-1,+1\}$ binary vector $x$ as in the Subset Sum problem. Then we have
\[ \frac{1-x_ix_j}{2} = \left\{ \begin{array}{rcl} 1,&\quad& i \in I,\ j \in \bar I\ \mbox{or}\ i \in \bar I,\ j \in I; \\
0,&& i,j \in I \mbox{or}\ i,j \in \bar I. \end{array} \right.
\]
The expression $\frac{1-x_ix_j}{2}$ can then be used as an indicator that the vertices $i,j$ lie in different subsets of the partition. Hence the cut weight can be written as
\[ \frac14 \sum_{i,j \in V} (1-x_ix_j)W_{ij} = \frac14(-x^TWx + {\bf 1}^TW{\bf 1}).
\]
Here the multiplier is $\frac14$ instead of $\frac12$ to compensate for the fact that every edge is counted twice in the double sum over $i,j$.

The MaxCut problem then can be reformulated equivalently as the QUBO problem
\[ \frac{1}{4} \min_{x \in \{-1,+1\}^n} \, x^T W x,
\]
up to a constant summand in the objective. Note that the sign disappeared because we converted a maximization problem to a minimization problem.

\medskip

\section{Methods for Solving QUBO} \label{sec:QUBOmethods}
\subsection{Metaheuristic Algorithms}

Metaheuristic algorithms are general-purpose strategies designed to explore large and complex optimization landscapes where exact methods are infeasible. They aim to balance exploration of the global search space with exploitation of promising regions. In this subsection, we discuss three representative metaheuristics: \emph{Simulated Annealing}, \emph{Tabu Search}, and the tensor-based \emph{PROTES} method.  

\subsubsection{Simulated Annealing (SA).}  
Simulated Annealing \cite{kirkpatrick1983optimization} is inspired by the physical process of annealing in metallurgy, where a material is heated and slowly cooled to achieve a low-energy crystalline state. Formally, given an objective (energy) function $E(x)$ defined over a discrete solution space $\mathcal{X}$, SA maintains a single candidate solution $x \in \mathcal{X}$ and performs stochastic updates. At iteration $t$, a neighbor $x'$ of the current state $x$ is sampled from a proposal distribution $\mathcal{N}(x)$, and the acceptance probability is defined by the Metropolis criterion:
\[
P(x \to x') = \min \left(1, \exp\left( -\frac{E(x') - E(x)}{T_t} \right) \right),
\]
where $T_t$ is a temperature parameter that decreases according to a cooling schedule, e.g.\ $T_t = \alpha^t T_0$ with $\alpha < 1$.  

\begin{algorithm}[H]
\caption{Simulated Annealing}
\KwIn{Initial solution $x_0$, temperature schedule $\{T_t\}$}
\For{$t = 0,1,2,\dots$ until stopping criterion}{
  Sample $x' \sim \mathcal{N}(x_t)$\;
  Compute $\Delta = E(x') - E(x_t)$\;
  \eIf{$\Delta \le 0$}{
    $x_{t+1} \gets x'$\;
  }{
    $x_{t+1} \gets x'$ with prob.\ $\exp(-\Delta/T_t)$, else $x_{t+1} \gets x_t$\;
  }
}
\KwOut{Best solution found}
\end{algorithm}

\subsubsection{Tabu Search (TS).}  
Tabu Search \cite{glover1989tabu} is a deterministic metaheuristic that enhances local search by using memory structures. At each iteration, the algorithm moves to the best admissible neighbor of the current solution, even if it leads to a worse objective value. To avoid cycles, recently visited solutions (or moves) are recorded in a \emph{tabu list}, which forbids revisiting them for a fixed number of iterations (the \emph{tabu tenure}). An aspiration criterion may override the tabu status if a move yields a solution better than any found so far.  

Formally, let $\mathcal{N}(x)$ be the neighborhood of $x$ and $\mathcal{T}_t$ the tabu list at iteration $t$. The next state is chosen as:
\[
x_{t+1} = \arg\min_{x' \in \mathcal{N}(x_t) \setminus \mathcal{T}_t} E(x').
\]

\begin{algorithm}[H]
\caption{Tabu Search}
\KwIn{Initial solution $x_0$, tabu tenure $L$}
Initialize tabu list $\mathcal{T} \gets \emptyset$\;
\For{$t = 0,1,2,\dots$ until stopping criterion}{
  Determine admissible neighbors $\mathcal{N}'(x_t) = \mathcal{N}(x_t) \setminus \mathcal{T}$\;
  Select $x_{t+1} \in \arg\min_{x' \in \mathcal{N}'(x_t)} E(x')$\;
  Update $\mathcal{T}$ with $x_t$ (expire entries older than $L$ steps)\;
}
\KwOut{Best solution found}
\end{algorithm}

\subsubsection{Probabilistic Optimization with Tensor Sampling.}  
The PROTES algorithm (Probabilistic Optimization with Tensor Sampling) \cite{batsheva2023protes} is a modern metaheuristic based on the tensor-train (TT) decomposition. Instead of maintaining a single trajectory (as in SA and TS), PROTES maintains a probability distribution $p(x)$ over $\mathcal{X} = \{0,1\}^d$ in TT format:
\begin{equation}    
p(x_1,\ldots,x_d) = G^{(1)}_{x_1} G^{(2)}_{x_2} \cdots G^{(d)}_{x_d},
\end{equation}

where each $G^{(k)}_{x_k}$ is a low-rank matrix slice from the $k$th TT-core. This representation enables efficient sampling and updating in exponential spaces.  

At each iteration, PROTES samples a batch of candidate solutions $\{x^{(j)}\}_{j=1}^M \sim p(x)$, evaluates their objective values $E(x^{(j)})$, and selects an elite subset $\mathcal{E}$ of the best samples. The TT-parameters are updated by maximizing the likelihood of $\mathcal{E}$, effectively shifting the distribution toward regions of lower energy:
\[
\theta \gets \arg\max_\theta \sum_{x \in \mathcal{E}} \log p_\theta(x).
\]

\begin{algorithm}[H]
\caption{High-Level TT-Sampling Workflow (PROTES)}
\label{alg:tt_sampling}
\begin{enumerate}
  \item \textbf{Initialization.}  
    Set up the TT sampler with prescribed TT–ranks.
  \item \textbf{Sampling.}  
    Draw \(K\) candidate points 
    \(\{x^{(j)}\}_{j=1}^K\)
    and evaluate \(f\bigl(x^{(j)}\bigr)\).
  \item \textbf{Selection.}  
    From these \(K\), pick the top \(k\) points with the highest function values.
  \item \textbf{Refinement.}  
    Perform a few steps of gradient ascent on the log-likelihood of the selected \(k\) samples,  
    using learning rate \(\lambda\).
\end{enumerate}
\end{algorithm}

Simulated Annealing explores the landscape stochastically with gradually decreasing randomness, Tabu Search uses deterministic moves with adaptive memory to prevent cycling, and PROTES employs a high-dimensional probabilistic model to capture global structures of the solution space. Together, these represent three different philosophies of metaheuristic design: stochastic annealing, memory-guided search, and tensor-based probabilistic optimization.  

\subsection{Quantum Annealing and Methods Inspired by it}

We begin with the quadratic unconstrained binary optimization (QUBO) problem:
\begin{align}
\underset{x \in \{0,1\}^n}{\text{min}} \ \sum_i a_i x_i + \sum_{i>j} b_{ij} x_i x_j,
\end{align}
where the vector $x$ encodes binary decision variables, $a_i$ are linear biases, and $b_{ij}$ represent pairwise couplings. As we mentioed before, QUBO formulations capture a wide range of NP-hard combinatorial optimization problems and naturally map onto the Ising Hamiltonian formalism.

\subsubsection{Quantum Annealing}

Quantum annealing (QA) is a metaheuristic designed to solve combinatorial optimization problems by harnessing controlled quantum dynamics. Conceptually, it generalizes classical simulated annealing (SA) \cite{ayanzadeh2020reinforcement}, but replaces thermal fluctuations with quantum fluctuations as the driver of state transitions. The method was introduced by Kadowaki and Nishimori \cite{kadowaki1998quantum}, who proposed using a transverse-field Ising model in which a gradually decreasing transverse field plays a role analogous to temperature in SA, but crucially enables quantum tunneling across energy barriers. This tunneling can allow QA to escape exponentially narrow local minima that would trap SA.

Modern QA is realized on superconducting hardware platforms such as D-Wave’s Advantage QPU, where arrays of flux qubits are arranged in sparse graphs with programmable couplers and biases. These devices implement an interpolation between two Hamiltonians: an initial “driver” Hamiltonian with a simple ground state and a final “problem” Hamiltonian encoding the QUBO instance. The quantum system is prepared in the easily accessible ground state of the driver and then evolved toward the problem Hamiltonian according to the adiabatic principle, ideally remaining close to the instantaneous ground state throughout.

\subsubsection{Workflow and Embedding}

In practice, QA proceeds in three phases: (i) embedding, (ii) annealing evolution, and (iii) readout.  
\begin{enumerate}
    \item \emph{Embedding.} The logical QUBO graph must be embedded onto the sparse hardware connectivity (e.g., the Pegasus topology in D-Wave systems). Logical variables are represented by chains of physical qubits linked by strong ferromagnetic couplings to enforce consistency. Embedding is itself an NP-hard problem and typically solved using heuristics. 
    \item \emph{Annealing Evolution.} The system evolves under a time-dependent Hamiltonian
    \begin{equation}
        H(s) = -\frac{A(s)}{2} \sum_i \sigma_x^{(i)} 
        + \frac{B(s)}{2} \left( \sum_i h_i \sigma_z^{(i)} 
        + \sum_{i>j} J_{ij} \sigma_z^{(i)} \sigma_z^{(j)} \right),
    \end{equation}
    where $s \in [0,1]$ parametrizes the annealing schedule, $A(s)$ decreases while $B(s)$ increases, and $\sigma_x, \sigma_z$ are Pauli operators. At $s=0$, the ground state is a uniform superposition dominated by the transverse field; at $s=1$, the Hamiltonian reduces to the classical Ising model representing the QUBO. 
    \item \emph{Readout.} After evolution, the system is measured in the computational basis, collapsing each qubit to a binary value. Repeated runs yield a distribution over low-energy solutions, from which the best solution is selected.
\end{enumerate}

\begin{figure}
\centering
\includegraphics[width=0.75\textwidth]{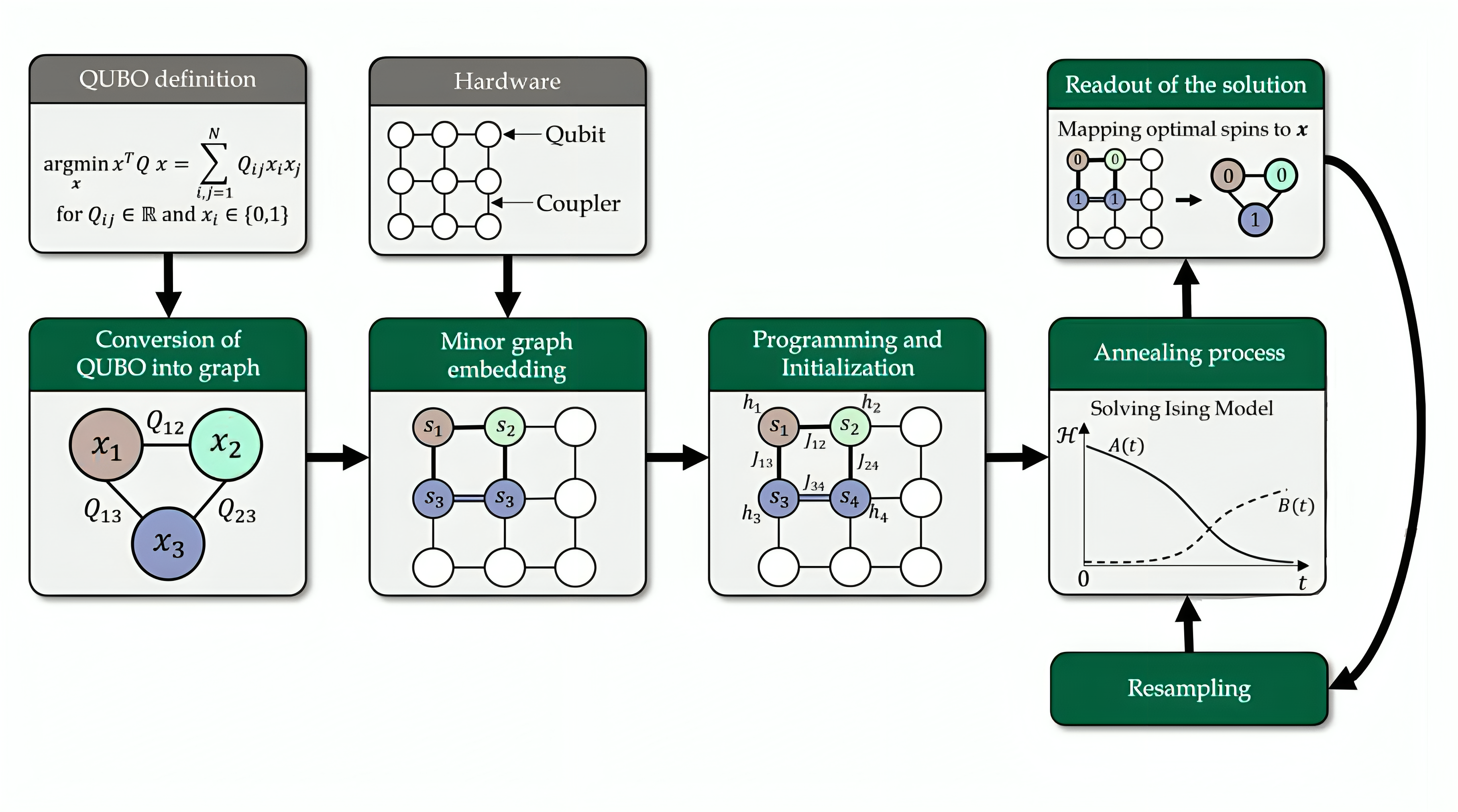}
\caption{
Workflow of quantum annealing for QUBO: 
(1) QUBO formulation $\text{argmin}\, x^T Q x$ and mapping to Ising form;
(2) Minor-embedding of the logical graph into the hardware graph (qubits and couplers);
(3) Programming of annealing schedule ($A(t)$: transverse field, $B(t)$: problem Hamiltonian, $J(t)$: coupler strengths);
(4) Quantum evolution under the time-dependent Hamiltonian;
(5) Measurement and decoding into binary solution $x$. Adapted from \cite{goodrich2018optimizing}.}
\label{fig:flowchart}
\end{figure}

This process is illustrated in Fig.~\ref{fig:flowchart}. In practice, non-adiabatic effects, thermal excitations, and control noise may prevent the system from reaching the true ground state, but QA often produces high-quality approximate solutions with polynomial overhead.

\subsubsection{Quantum-Inspired Classical Methods}

Although QA offers an elegant physical paradigm, it requires cryogenic quantum hardware with limited qubit counts and restricted connectivity. To broaden applicability, researchers have developed \emph{quantum-inspired annealing} techniques that emulate the essential physics of QA using classical algorithms. Two notable approaches are \emph{Simulated Quantum Annealing} (SQA) and \emph{Symmetry Breaking} (SB).
SQA reproduces the dynamics of QA via classical stochastic simulation, most prominently using path-integral Monte Carlo (PIMC) methods \cite{kadowaki1998quantum,morita2008mathematical}. 

The system is evolved under an interpolating Hamiltonian
\begin{equation}
    H_{\mathrm{SQA}}(t) = A(t) H_0 + B(t) H_P,\quad t\in[0,T],
\end{equation}
with $A(0)=1, B(0)=0$ and $A(T)=0, B(T)=1$. The quantum partition function is mapped to a classical effective model via a Suzuki–Trotter decomposition, introducing an additional “imaginary-time” dimension of $M$ Trotter slices. Each slice corresponds to a classical replica of the system, coupled along the imaginary-time axis to emulate quantum fluctuations. Monte Carlo updates over this extended ensemble allow the system to undergo \emph{quantum-like tunneling} between energy minima, thereby escaping barriers more efficiently than thermal hopping alone. The algorithm’s accuracy depends on the annealing schedule $(A(t),B(t))$, the number of Trotter slices $M$, and the choice of Monte Carlo dynamics. 

Simulated Bifurcation (SB) takes a deterministic, dynamical-systems perspective. Instead of stochastic sampling, it encodes the QUBO in a system of coupled differential equations:
\begin{equation}    
\dot{x}_i = y_i,\quad
\dot{y}_i = -(\mu(t)-\lambda)x_i + \sum_j J_{ij}x_j - \gamma y_i,
\end{equation}
where $x_i$ are continuous state variables that gradually bifurcate toward $\pm 1$ values. As the control parameter $\mu(t)$ ramps through the critical threshold $\mu_c = \lambda + \rho(J)$ (with $\rho(J)$ the spectral radius of $J$), the system undergoes a pitchfork bifurcation that spontaneously breaks symmetry, producing binary states. The quality of solutions depends sensitively on the damping parameter $\gamma$, offset $\lambda$, and the annealing schedule for $\mu(t)$. SB thus replaces stochastic tunneling with nonlinear dynamical instabilities as the exploration mechanism.

SQA and SB embody two complementary interpretations of annealing-inspired optimization:
\begin{itemize}
  \item \textbf{SQA} leverages \emph{stochastic sampling} in an extended configuration space to mimic quantum tunneling.
  \item \textbf{SB} exploits \emph{deterministic bifurcations} in nonlinear dynamics to drive continuous states toward binary minima.
\end{itemize}
Both can be abstracted as probabilistic oracles producing near-optimal solutions with high probability:
\begin{equation}
\mathbb{P}\!\left(H_{\mathrm{QUBO}}(\mathbf{x}^*) \le \min_{\mathbf{x}} H_{\mathrm{QUBO}}(\mathbf{x}) + \delta\right) \ge 1 - \eta,
\end{equation}
where $\mathbf{x}^*$ is the returned configuration, $\delta$ the approximation gap, and $\eta$ the failure probability.

Quantum annealing introduces a physically motivated optimization framework that leverages quantum tunneling and adiabatic evolution. Its quantum-inspired counterparts, SQA and SB, transplant these ideas into classical settings—SQA by simulating quantum statistical mechanics, and SB by exploiting dynamical instabilities. Together, they extend the annealing paradigm beyond hardware constraints, offering scalable methods for tackling large-scale QUBO problems while illuminating the interplay between physics and optimization.

\subsection{Burer's rank 2 relaxation}

In this section we present another method which inspired the present work, namely Burer's rank 2 relaxation \cite{Burer02}. It relaxes the rank 1 constraint on the semi-definite matrix $X$ to a rank 2 constraint and solves the non-convex problem
\begin{equation}\label{burer}
    \min_X\,\langle Q, X\rangle: \qquad X \in \mathcal{S}_+^n, \diag\,X = 1,\ \rk(X) \leq 2.
\end{equation}
The key idea underlying the relaxation is that any feasible matrix in this problem can be written in the form
\[
X = (\cos(\varphi_i-\varphi_j))_{i,j=1,\dots,n}
\]
for some angles $\varphi_1,\dots,\varphi_n$. The objective then becomes
\[ \langle Q, X\rangle = \sum_{i,j=1}^n Q_{ij}\cos(\varphi_i-\varphi_j),
\]
i.e., it reduces to a trigonometric polynomial on the $n$-dimensional torus. This problem is non-convex but smooth, in addition the derivatives of the objective of arbitrary order can easily be computed. The minimization can hence be performed by iterative algorithms like gradient descent or the Newton method.

The method, of course, does not guarantee convergence to the global minimum. In general, we obtain a local minimum $\varphi^*$ of the objective. The corresponding matrix $X^*$ can then be used as a covariance matrix in the GW randomized procedure. 

Note that the rows of the factor $F$ in the factorization $X = FF^T$ are of the form $(\cos\phi_i,\sin\phi_i)$, and hence the geometric interpretation of the GW procedure just partitions given $n$ points on the unit circle by an isotropic line through the origin. All such partitions can be enumerated in linear time and the best can be chosen.

In our method we adopt the idea of minimizing the objective over the factor $F$ as independent variable rather than the positive semi-definite matrix $X$.

\subsubsection{Special case $k = 2$}

If the rank of $X$ does not exceed 2, as in the relaxation proposed by Burer et al., the cost function \eqref{expectationFormula2} becomes piece-wise affine.

In this case we have $X_{ij} = \cos(\varphi_i - \varphi_j)$ for some angles $\varphi_1,\dots,\varphi_n$. The optimization can be accomplished directly in the space of these angles. Let $\delta \in \mathbb R^n$ be a direction in angle space. Let us compute the cost for the updated angles $\varphi + \delta$. 

Let $\tau_{ij} = \frac{2}{\pi}\arcsin X_{ij} \in [-1,1]$. Then 
\[ \tau_{ij} = 1 + 4l + \sigma\frac{2}{\pi}(\varphi_i - \varphi_j),
\]
where the sign $\sigma \in \{-1,1\}$ and the number $l \in \mathbb Z$ are chosen such that the right-hand side is in the interval $[-1,1]$ (note that changing an angle by a multiple of $2\pi$ changes $l$ by an integer). Then we get for small enough $\delta$
\[ \tau_{ij}(\varphi + \delta) = \tau_{ij}(\varphi) + \left\{ \begin{array}{rcl} \sigma\frac{2}{\pi}(\delta_i - \delta_j),&\quad& \tau_{ij} \in (-1,1), \\
- \frac{2}{\pi}|\delta_i - \delta_j|,&& \tau_{ij} = 1, \\ \frac{2}{\pi}|\delta_i - \delta_j|,&& \tau_{ij} = -1. \end{array} \right.
\]

It follows that the directional derivative is piece-wise linear in $\delta$. It can also be minimized by a difference of convex functions minimization algorithm. Again the difficult step is the update of $x_k$ by an auxiliary minimization problem. In this case, however, this step amounts to solving a linear program. 

\section{Proposed approach}

We present a new method for solving Quadratic Unconstrained Binary Optimization (QUBO). We look for a covariance matrix with unit diagonal that directly minimizes the expected value of the solution obtained by Goemans-Williamson randomized rounding using this covariance matrix. This is in contrast to the semi-definite relaxation, where the objective function itself, or rather its linear extension from the set of rank 1 matrices to the whole spectrahedron ${\cal SR}$ is minimized. This alignment between the continuous surrogate and the discrete target yields both algorithmic and theoretical advantages that we now make explicit.

Let $Q\in\mathbb{R}^{n\times n}$ be the symmetric matrix defining the QUBO objective and let $X$ denote a feasible solution of the canonical SDP relaxation (so $\operatorname{diag}(X)=\mathbf{1}$ and $X\succeq 0$). The Goemans–Williamson rounding scheme implies a closed-form expression for the expected objective value produced by a Gaussian sign rounding \cite{Nesterov98}:

$$
\mathbb{E}_\xi\big[\,\langle Q, \operatorname{sgn}(\xi^\top F)\operatorname{sgn}(\xi^\top F)^\top\rangle\,\big]
 \;=\; \frac{2}{\pi}\,\langle Q,\arcsin(X)\rangle,\qquad X = FF^\top,
$$

where $\arcsin$ is applied element-wise and $F\in\mathbb{R}^{n\times k}$ is any factor satisfying $FF^\top=X$. 

Our methods' main idea is to directly minimize the nonconvex objective

\begin{equation} \label{expectationFormula2}
\Phi(F)\;=\;\frac{2}{\pi}\,\langle Q,\arcsin(FF^\top)\rangle,    
\end{equation}
which measures the expected post-rounding performance of the factor $F$. Optimizing $\Phi$ over the factorization manifold $\{F:\,\|F_{i:}\|_2=1\ \forall i\}$ has two important consequences: (i) the optimization is low-rank (controlled by $k$) and thus scalable in practice when $k\ll n$; (ii) progress in $\Phi$ has an immediate interpretability in terms of the combinatorial objective after rounding. Next we detail how to perform the minimization of $\Phi$.

\subsection{Computing the directional derivative}

We consider the problem of minimization of expectation \eqref{expectationFormula2} over factors $F \in \mathbb R^{n \times k}$ of matrices $X = FF^T \succeq 0$ with unit diagonal not exceeding rank $k$. Then the problem becomes
\[ \max_{F \in \mathbb R^{n \times k}} \frac{2}{\pi}\langle Q,\arcsin (FF^T) \rangle: \qquad \sum_{j=1}^k F_{ij}^2 = 1 \quad \forall\ i = 1,\dots,n.
\]

The gradient of the cost function is given by
\[ G = \frac{4}{\pi}\frac{Q}{\sqrt{1-(FF^T)^2}}F,
\]
where in the expression $\frac{Q}{\sqrt{1-(FF^T)^2}}$ the operations are meant element-wise, i.e.,
\[ \left( \frac{Q}{\sqrt{1-(FF^T)^2}} \right)_{ij} = \frac{Q_{ij}}{\sqrt{1-\langle F_{i*},F_{j*} \rangle^2}}.
\]
The gradient will be restricted to the subspace given by $\diag(FG^T) = 0$ which is tangent to the variety of matrices $F$ which satisfy the quadratic constraints on the row norms of $F$.

Note that the gradient is not defined when an off-diagonal element of $X = FF^T$ is equal to $\pm1$.  Consider such an element $X_{ij} = \sigma \in \{-1,+1\}$ in more detail. Let $D \in \mathbb R^{n \times k}$ be a direction, and suppose that two rows $f_i,f_j$ of $F$ obey $\langle f_i,f_j \rangle = \sigma$. As a consequence $f_j = \sigma f_i$. Let $d_i,d_j$ be the corresponding rows of $D$. Then $\langle f_i,d_i \rangle = \langle f_i,d_j \rangle = 0$. The normalized rows of the updated factor $F + tD$, coming from a step in the direction $D$ with step size $t$, where we suppose $t \geq 0$, are given by
\[ \frac{f_i + td_i}{\sqrt{\langle f_i + td_i,f_i + td_i \rangle}}, \quad \frac{\sigma f_i + td_j}{\sqrt{\langle \sigma f_i + td_j,\sigma f_i + td_j \rangle}}.
\]
Their scalar product amounts to
\begin{equation} \label{XijExpression}
    X_{ij}(t) = \frac{\sigma + t^2\langle d_i,d_j \rangle}{\sqrt{1 + t^2\langle d_i,d_i \rangle}\sqrt{1 + t^2\langle d_j,d_j \rangle}} = \sigma + t^2\langle d_i,d_j \rangle - \frac{t^2}{2}\sigma\langle d_i,d_i \rangle - \frac{t^2}{2}\sigma\langle d_j,d_j \rangle + O(t^4).
\end{equation} 

Let us compute the directional derivative
\[ \mu = \lim_{t \to +0}\frac{\frac{2}{\pi}\arcsin X_{ij}(t) - \frac{2}{\pi}\arcsin X_{ij}(0)}{t}
\]
of the expression $\frac{2}{\pi}\arcsin X_{ij}$ in the direction $D$. By definition we get $\frac{2}{\pi}\arcsin X_{ij}(t) = \sigma + t\mu + O(t^2)$, where $\sigma\mu \leq 0$, and
\[ X_{ij}(t) = \sigma - \frac{\pi^2t^2\sigma}{8}\mu^2 + O(t^4).
\]
By comparison with \eqref{XijExpression} it follows that
\[ \mu = -\frac{2}{\pi}\sigma\|d_i - \sigma d_j\|.
\]
In particular, the directional derivative of $\frac{2}{\pi}\arcsin X_{ij}$ in the direction $D$ equals that in the direction $-D$ instead of being opposite as in the case of $C^1$ functions. However, for every direction $D$ the directional derivative of the cost function exists and is given by
\[ G_D = 2\sum_{i < j}Q_{ij}\mu_{ij}, \qquad \mu_{ij} = \left\{ \begin{array}{rcl} \frac{2}{\pi}\frac{1}{\sqrt{1-X_{ij}^2}}(e_i^T(FD^T+DF^T)e_j),&\quad& X_{ij} \in (-1,+1), \\ \frac{2}{\pi}\|d_i + d_j\|,&& X_{ij} = -1, \\ 
-\frac{2}{\pi}\|d_i - d_j\|,&& X_{ij} = 1. \end{array} \right.
\]
Here $e_i$ are the standard basis vectors in $\mathbb R^n$.

Finally we arrive at
\begin{align}
\nabla_D \Phi(F) \;=\; \left.\frac{\mathrm{d}}{\mathrm{d}\varepsilon}\right|_{\varepsilon=0}\Phi(F+\varepsilon D)
 \;= 
\begin{cases}
\displaystyle \frac{2}{\pi}\sum_{i<j} Q_{ij}\,
    \frac{f_{i}d_{j}^\top + f_{j}d_{i}^\top}{\sqrt{1 - X_{ij}^2}}, 
    & \text{if } |X_{ij}| < 1, \\[2.0ex]
\displaystyle \frac{2\sigma}{\pi}\sum_{i<j} Q_{ij}\,
    \|d_{i} + \sigma d_{j}\|, 
    & \text{if } X_{ij} = \sigma \in \{1,-1\}.
\end{cases}
\end{align}

Note that in case $X_{ij} \in \{-1,+1\}$, if $Q_{ij}X_{ij} > 0$, then the function $\frac{2}{\pi}Q_{ij}\arcsin X_{ij}$ has a maximum, if $Q_{ij}X_{ij} < 0$, then it has a minimum. 

In order to apply a gradient method to minimize cost function \eqref{expectationFormula2}, we need to find a descent direction in order to make a step. The directional derivative, as a function of the direction $D$, consists of three parts:
\begin{itemize}
  \item a linear part $2\sum_{i < j: |X_{ij}| < 1}Q_{ij}\mu_{ij}$
  \item a non-positive concave conic part $2\sum_{i < j: X_{ij} = \operatorname{sgn} Q_{ij}}Q_{ij}\mu_{ij}$
  \item a nonnegative convex conic part $2\sum_{i < j: X_{ij} = -\operatorname{sgn} Q_{ij}}Q_{ij}\mu_{ij}$
\end{itemize}
Minimizing the directional derivative over $D$ can be accomplished by an algorithm used for minimization of differences of convex functions, described in Section \ref{sec:diffConvex}. Here the concave part is approximated by a linear upper bound, while the linear and convex part form a convex sum.

The difficult step in this scheme is hence the search for a descent direction in order to update $x_k$, which can be accomplished by solving an auxiliary convex minimization problem. This problem amounts to minimizing the sum of a linear function and some norms of linear combinations of rows of $D$. This can be written as a second-order cone program (SOCP).


\subsection{Minimization of differences of convex functions} \label{sec:diffConvex}

We briefly outline the idea of difference-of-convex (DC) minimization. Consider the problem
\[ \min_x (g(x) - h(x))
\]
where $g,h$ are convex functions. This is a non-convex optimization problem. 

Consider the following algorithm \ref{alg:DCminimization}.

\medskip

\begin{algorithm}[H]
\caption{Difference of convex functions minimization}
\label{alg:DCminimization}
\KwIn{\(x_0,N\)}
\KwOut{\(x_N\)}
\For{\(k=1:N\)}{
  \(y_k \gets \nabla h(x_{k-1})\)
  \(x_k \gets \arg\min_x(g(x) - \langle y_k,x \rangle)\)
}
\Return \(x_N\)
\end{algorithm}

\medskip

Here $y_k$ is in the subgradient of $h$ at $x_{k-1}$.

Since $h$ is convex, we have
\[ h(x_k) \geq h(x_{k-1}) + \langle y_k,x_k-x_{k-1} \rangle.
\]
Further, we have for every $x$ that
\[ g(x_k) - \langle y_k,x_k \rangle \leq g(x) - \langle y_k,x \rangle.
\]
In particular,
\[ g(x_k) - h(x_k) \leq g(x_{k-1}) + \langle y_k,x_k-x_{k-1} \rangle - h(x_k) \leq g(x_{k-1}) - h(x_{k-1}),
\]
and the sequence of function values is monotonically decreasing. 

Clearly local minima of the cost function are stationary points of the algorithm.

\subsection{Direct Expectation Minimization}
Algorithm \ref{alg:SOCP} below describes the computation of a descent direction $D$ for the Direct Expectation Minimization (DEM) update.  The direction is obtained by solving a second-order cone program (SOCP) that exactly captures the three geometric contributions that appear in the directional derivative of the expected objective.  Pairs $(i,j)$ with $X_{ij}= (FF^\top)_{ij}\in(-1,1)$ produce a linear contribution in $D$.  Pairs on the manifold boundary $|X_{ij}|\approx 1$ produce non-linear, nonsmooth parts: when the boundary part is convex ($X_{ij} = -\operatorname{sgn}(Q_{ij})$) it is represented by a second order cone (SOC) term; concave boundary parts ($X_{ij} = \operatorname{sgn}(Q_{ij})$) are treated within the SOCP via a subgradient linearization (using information from the previous direction).  The SOCP therefore combines the assembled linear term, SOC constraints for convex boundary parts, and subgradient constraints for concave boundary parts, subject to the tangency constraints that keep the update on the manifold and a norm bound that fixes the descent scale.

\begin{algorithm}[t]
\caption{Compute descent direction \(D\) (SOCP auxiliary problem)}
\label{alg:SOCP}
\KwIn{\(F\in\mathbb{R}^{n\times k}\) (rows \(f_i^\top\)), \(Q\in\mathbb{S}^n\), \(D^{\mathrm{curr}}\in\mathbb{R}^{n\times k}\)}
\KwOut{Solution \(D\in\mathbb{R}^{n\times k}\)}
\BlankLine
\(X \leftarrow F F^\top\) \\
Introduce epigraph variables \(v_{ij}\) for all \(i<j\) \\
Build linear objective: \(\displaystyle \min_{D,\;v}\; \sum_{i<j} v_{ij}\) \\
Add orthogonality constraints \(\langle f_i, d_i\rangle = 0\) for all \(i\) \\
Add norm bound \(\|D\|_F \le 1\) (Frobenius norm) \\
\For{each pair \(i<j\)}{
  \(\alpha_{ij}\leftarrow \tfrac{2}{\pi} Q_{ij}\) \\
  \eIf{\(|X_{ij}|<1-\varepsilon\)}{ 
    Add linear equality \\
    \(\displaystyle v_{ij} = \frac{\alpha_{ij}}{\sqrt{1-X_{ij}^2}}\big(f_i^\top d_j + f_j^\top d_i\big)\) \Comment*[r]{interior (smooth) case}
  }{
    \eIf{\(|X_{ij}-(-1)|<\varepsilon\)}{ 
      \eIf{\(X_{ij} = -\operatorname{sgn}(Q_{ij})\)}{ 
        Add second-order cone constraint \\
        \(v_{ij} \;\ge\; \alpha_{ij}\,\|\,d_i + d_j\,\|_2\) \Comment*[r]{convex boundary}
      }{
        Compute \(s_{ij}\in\partial\|d_i^{\mathrm{curr}}+d_j^{\mathrm{curr}}\|_2\) \\
        Add affine constraint \\
        \(v_{ij} \;\ge\; \alpha_{ij}\,\langle s_{ij},\,d_i + d_j\rangle\) \Comment*[r]{concave boundary}
      }
    }{
      \If{\(|X_{ij}-1|<\varepsilon\)}{ 
        \eIf{\(X_{ij} = -\operatorname{sgn}(Q_{ij})\)}{ 
          Add second-order cone constraint \\
          \(v_{ij} \;\ge\; \alpha_{ij}\,\|\,d_i - d_j\,\|_2\) \Comment*[r]{convex boundary}
        }{
          Compute \(s_{ij}\in\partial\|d_i^{\mathrm{curr}}-d_j^{\mathrm{curr}}\|_2\) \\
          Add affine constraint \\
          \(v_{ij} \;\ge\; \alpha_{ij}\,\langle s_{ij},\,d_i - d_j\rangle\) \Comment*[r]{concave boundary}
        }
      }
    }
  }
}
Solve the resulting convex SOCP for \((D,v)\) \\
\Return \(D\)
\end{algorithm}

\subsubsection{Exact DEM}
The exact DEM approach uses a projected gradient method, where a second-order cone programming (SOCP) auxiliary problem is used to figure out the direction of each step.  Starting from an initial iterate $F^{(0)}$, the algorithm first computes a direction $D$ by solving the SOCP, ensuring that $D$ satisfies both the tangency and unit-norm constraints imposed by the auxiliary problem. The iterate is then updated explicitly according to $F \leftarrow F + \eta D$, where $\eta$ denotes the step size. Following this update, $F$ is retracted onto the feasible manifold by normalizing each row individually, that is, $f_i \leftarrow f_i / \|f_i\|_2$ for all rows $i$, thereby maintaining feasibility with respect to the manifold constraints. The iterative procedure continues until the computed search direction $D$ becomes (nearly) zero in norm, or the maximum number of iterations is reached. The initial direction can be chosen randomly and normalized, and optional periodic evaluation or printing of the objective function may be performed to monitor convergence.

\begin{algorithm}[H]
\caption{Exact DEM}
\label{alg:exact-DEM}
\KwIn{Initial $F^{(0)}$, $Q$, tolerance $\text{tol}$, step size $\eta$, max iterations $T$}
\KwOut{Optimized $F$}
$F \gets F^{(0)}$ \\
Initialize $D \gets$ random matrix in $\mathbb{R}^{n \times k}$, then normalize $D \gets D / \|D\|_F$ \\
\For{iteration $t = 1$ \KwTo $T$}{
    Compute descent direction $D \gets \text{solve\_socp}(F, Q, D)$ (Algorithm~\ref{alg:SOCP}) 
    \BlankLine
    \Comment*[r]{Returns $D$ (reshaped $n \times k$)}
    Take explicit step: $F \gets F + \eta D$\\
    Row-wise normalization: $f_i \gets \frac{f_i}{\|f_i\|_2}$ for $i = 1, \dots, n$\\
    Evaluate objective $\Phi(F)$\\
    \If{$\|D\|_F < \text{tol}$}{
        break \Comment*[r]{Direction vanishes — stop}
    }
}
\Return \(F\)
\end{algorithm}



\subsection{Direct Expectation Minimization with Riemannian Descent and Clipping (DEM-RC)}

The \emph{Direct Expectation Minimization} (DEM) approach approximates this discrete optimization problem by operating in a continuous embedding space. In particular, we represent binary variables through low-rank spherical embeddings $F \in \mathbb{R}^{n \times r}$ with row normalization. The expected objective can be expressed as a smooth function of $X = FF^\top$, and optimized via Riemannian gradient descent. 

To ensure numerical stability, the entries of $X$ are clipped to the interval $(-1+\varepsilon,\, 1-\varepsilon)$. This prevents singularities when evaluating the gradient
\[
    G \;=\; \frac{2}{\pi} \left( Q \oslash \sqrt{1 - X^{\circ 2}} \right) F,
\]
where $\oslash$ denotes elementwise division and $X^{\circ 2}$ the Hadamard (element-wise) square. After each update step
\[
    F \;\gets\; \text{Normalize}\!\left(F - \eta G \right),
\]
row-wise normalization maintains the embedding on the unit sphere.

Finally, to extract a feasible binary solution, we apply a Goemans--Williamson-style randomized rounding: multiple Gaussian vectors $g \sim \mathcal{N}(0,I_r)$ are sampled, and binary candidates are generated via $x = \mathrm{sign}(Fg)$. The best solution across $R$ trials is returned.

The procedure is summarized in Algorithm~\ref{alg:demrc}.

\begin{algorithm}[H]
\caption{DEM-RC: Direct Expectation Minimization with Riemannian descent and Clipping}
\label{alg:demrc}
\KwIn{QUBO matrix $Q \in \mathbb{R}^{n \times n}$, rank $r$, steps $T$, rounds $R$, step size $\eta$, clipping $\varepsilon$}
\KwOut{Binary solution $x^* \in \{-1,1\}^n$}

\BlankLine
\textbf{Initialization:} \\
\For{$i=1$ \KwTo $n$}{
    Sample $f_i \sim \mathcal{N}(0,I_r)$ \\
    Normalize $f_i \gets f_i / \|f_i\|_2$
}
Form $F \gets [f_1; f_2; \ldots; f_n] \in \mathbb{R}^{n \times r}$ \\
Set best objective $E^* \gets +\infty$, best solution $x^* \gets \emptyset$

\BlankLine
\textbf{Riemannian gradient descent:} \\
\For{$t = 1$ \KwTo $T$}{
    $X \gets F F^\top$ \\
    Clip $X_{ij} \in (-1+\varepsilon,\,1-\varepsilon)$ for all $i,j$ \\

    $D \gets \sqrt{\,1 - X^{\circ 2}\,}$ \hfill (elementwise) \\
    $W \gets Q \oslash D$ \hfill (Hadamard division) \\
    $G_e \gets \tfrac{2}{\pi}\, W F$ \hfill (Euclidean gradient) \\

    \For{$i=1$ \KwTo $n$}{
        $s_i \gets G_{e,i}^\top f_i$ \\
        $g_i \gets G_{e,i} - s_i f_i$ \hfill (Riemannian projection) \\
        $f_i \gets f_i - \eta g_i$ \hfill (gradient step) \\
        $f_i \gets f_i / \|f_i\|_2$ \hfill (retraction to sphere)
    }
}

\BlankLine
\textbf{Randomized rounding (Goemans--Williamson style):} \\
\For{$k=1$ \KwTo $R$}{
    Sample $g \sim \mathcal{N}(0,I_r)$ \\
    $y \gets F g$ \\
    $x \gets \mathrm{sign}(y)$ \\
    $E \gets x^\top Q x$ \\
    \If{$E < E^*$}{
        $E^* \gets E$, \quad $x^* \gets x$
    }
}

\Return $x^*$
\end{algorithm}


\section{Results}
We evaluate our approach on two representative families of QUBO instances, chosen to capture both synthetic and structured problem classes:
\begin{enumerate}
    \item \textbf{Random dense QUBOs:} Symmetric matrices \(Q \in \mathbb{R}^{n \times n}\) with entries sampled i.i.d.\ from \(\mathcal{N}(0,1)\). We report results for problem sizes \(n = 200, 500, 1000\). These instances are unstructured but widely used as stress tests for scalability and robustness. 
    \item \textbf{IK problems:} Inverse kinematics formulations arising from robotic arm positioning, which lead naturally to QUBO formulations after discretization (\cite{salloum2025quantum}). These instances test whether algorithms can handle structured, application-driven energy landscapes.
\end{enumerate}

Tables~\ref{tab:random_qubo_sizes} and~\ref{tab:IK_results} summarize performance across methods. We compare classical relaxations (GW-SDP, Burer-2), metaheuristics (SA, SQA, Tabu Search, Simulated Bifurcation), and our DEM variants. For each, we report the best objective value attained and the corresponding runtime. Several key trends can be observed from the results:
\begin{itemize}
    \item \textbf{Relaxations.} The relaxation methods (GW-SDP and Burer-2) yield feasible but weaker solutions than heuristics, and their runtime grows poorly with $n$, especially GW-SDP due to the high cost of semidefinite programming. semidefinite programming.
    \item \textbf{Heuristics.} Classical metaheuristics such as Simulated Annealing (SA) and Tabu Search reliably achieve the best or near-best objective values across all tested sizes, surpassing relaxation methods both in terms of solution quality and computational efficiency. Simulated Bifuraction demonstrates competitive performance with consistently low runtimes, indicating strong scalability. In contrast, Simulated Quantum Annealing (SQA) attains high-quality solutions but at the expense of significantly higher runtime, limiting its practicality for large-scale instances.
    \item \textbf{DEM Variants.} The proposed DEM family exhibits a favorable quality–efficiency trade-off. In particular, DEM-RC-R10 provides the most balanced performance, achieving solution quality comparable to or exceeding that of heuristic baselines, while requiring only a fraction of their runtime. Low-rank clipped variants (R5, R10) demonstrate extremely fast convergence for moderate sizes (\(n=200\)) and remain competitive as the problem dimension increases.
\end{itemize}



\begin{table*}[ht]
  \centering
  \setlength{\tabcolsep}{6pt} 
  \renewcommand{\arraystretch}{1.2} 
  \caption{\textbf{Random QUBO}. Performance on synthetic instances with Gaussian random matrices. Reported are the best objective and runtime (in seconds). Best results are in \textbf{bold}, second-best are \underline{underlined}.
  Some entries are omitted because the corresponding runs exceeded practical time limits.}
  \label{tab:random_qubo_sizes}
  \begin{tabular}{l
                  rr
                  rr
                  rr}
    \toprule
    & \multicolumn{2}{c}{$n=200$} & \multicolumn{2}{c}{$n=500$} & \multicolumn{2}{c}{$n=1000$} \\
    \cmidrule(lr){2-3}\cmidrule(lr){4-5}\cmidrule(lr){6-7}
    Method & Objective & Time
           & Objective & Time
           & Objective & Time \\
    \midrule
    \multicolumn{7}{l}{\textbf{Relaxations}} \\
    GW-SDP        & -2878.45 & 3.02 
                  & -11047.79 & 64.30 
                  & -30842.91 & 389.11 \\
    Burer-2       & -2812.66 & 2.54  
                  & -11375.48 & 115.61  
                  & -32158.87 & 321.86  \\
    \addlinespace
    \midrule
    \multicolumn{7}{l}{\textbf{Heuristics}} \\
    SA            & \textbf{-2977.62} & 4.42   
                  & \textbf{-11848.06} & 32.56  
                  & \underline{-33910.53} & 80.61   \\
    SQA           & \underline{-2958.99} & 130.44   
                  & \underline{-11603.67} & 838.84  
                  & -- & --   \\
    Tabu          & \textbf{-2977.62} & 21.07   
                  & \textbf{-11848.06} & 24.82  
                  & \textbf{-33956.17} & \underline{29.14}      \\
    SB         & -2957.19 & 16.53   
                  & -11396.44 & \textbf{10.47}  
                  & -32500.21 & \textbf{14.02}   \\
    \addlinespace
    \midrule
    \multicolumn{7}{l}{\textbf{DEM Variants}} \\
    DEM-RC-R5     & -2899.85 & \underline{0.33}
                  & -11289.44 & 21.65
                  & -32565.58 & 61.41 \\
    DEM-RC-R10    & -2928.86 & \textbf{0.26}
                  & -11440.78 & \underline{16.94}
                  & -32978.01 & 61.71 \\
    DEM-R5        & -2852.81 & 606.11
                  & -11372.75 & 183.22
                  & -32073.72 & 180.14 \\

    DEM-R10       & -2921.91 & 30.65
                  & -11479.73 & 360.06
                  & -32795.11 & 624.66 \\
    \bottomrule
  \end{tabular}
\end{table*}

\begin{table*}[ht]
  \centering
  \setlength{\tabcolsep}{6pt} 
  \renewcommand{\arraystretch}{1.2} 
  \caption{\textbf{Inverse kinematic}. Comparison of methods for solving IK problems of different sizes. Best results are in \textbf{bold}, second-best are \underline{underlined}. Some entries are omitted because the corresponding runs exceeded practical time limits.
  }
  \label{tab:IK_results}
  \begin{tabular}{l
                  rr
                  rr
                  rr}
    \toprule
    & \multicolumn{2}{c}{$n=50$} & \multicolumn{2}{c}{$n=100$} & \multicolumn{2}{c}{$n=200$} \\
    \cmidrule(lr){2-3}\cmidrule(lr){4-5}\cmidrule(lr){6-7}
    Method & Objective & Time
           & Objective & Time
           & Objective & Time \\
    \midrule
    \multicolumn{7}{l}{\textbf{Relaxations}} \\
    GW-SDP        & \textbf{-1227.61} & 390.79 
                  & - & -
                  & - & - \\
    Burer-2       & -1163.71 & 5.18  
                  & -1078.37 & 144.05
                  & -997.78 & 858.35  \\
    \addlinespace
    \midrule
    \multicolumn{7}{l}{\textbf{Heuristics}} \\
    SA            & \textbf{-1227.61} & \textbf{2.2} 
                  & \textbf{-1188.14} & \textbf{2.29}  
                  & \underline{-1006.8} & \textbf{5.73}   \\
    SQA           & \textbf{-1227.61} & 21.88   
                  & \textbf{-1188.14} & 67.58  
                  & -1006.63 & 256.07  \\
    Tabu          & \textbf{-1227.61} & 21.02   
                  & \textbf{-1188.14} & 21.07  
                  & \textbf{-1007.01} & 21.14      \\
    SB         & -1223.89 & \underline{6.32}   
                  & -1184.78 & \underline{6.68}  
                  & \textbf{-1007.01} & \underline{6.65}   \\
    \addlinespace
    \midrule
    \multicolumn{7}{l}{\textbf{DEM Variants}} \\
    DEM-RC-R20     & \underline{-1224.22} & 9.06
                  & \underline{-1187.7} & 16.01
                  & -1005.84 & 100.74 \\
    DEM-R20        & \underline{-1224.22} & 129.59
                  & \underline{-1187.7} & 128.67 
                  & -1005.78 & 126.45 \\

    \bottomrule
  \end{tabular}
\end{table*}



\subsection{Relaxation-Based Mathematical Programming and Approximation Algorithms} \label{RelaxationBasedMethods}

Formally a QUBO problem is a mixed-binary quadratic program. As such it can be solved by general purpose branch and bound methods using convex quadratic relaxations and branching on binary variables to build a search tree. Such methods, however, disregard the special structure of QUBO, which lacks constraints other than the binary constraint on the variables.

We shall consider a more sophisticated relaxation in some detail, namely the semi-definite relaxation. It can be obtained in the following way.

Consider a QUBO problem with $\{-1,+1\}$ binary variables,
\[
\min_{x \in \{-1,1\}^n}\,x^TQx = \min_{x \in \{-1,1\}^n}\,\langle Q,xx^T \rangle.
\]
Let us define the \emph{MaxCut polytope}
\[ 
{\mathcal MC} = \conv\{ X = X^T \in \{-1,1\}^{n \times n} \mid X \succeq 0,\ \rk\,X = 1 \}.
\]
It has $2^{n-1}$ vertices, each of the form $X = xx^T$ with $x \in \{-1,1\}^n$. 

Then the QUBO problem can be written as the formally linear program
\[
 \max_{X \in {\mathcal MC}}\,\langle Q,X \rangle.   
\]

The \textit{semi-definite relaxation} replaces the exponentially complex polytope ${\mathcal MC}$ by the simpler spectrahedron
\[ {\mathcal SR} = \{ X \succeq 0 \mid \diag(X) = \textbf{1} \} \supset {\mathcal MC}
\]
which constitutes an outer approximation of ${\mathcal MC}$. Minimizing the objective over this set is equivalent to the semi-definite program
\[ \max_{X \succeq 0} \langle Q,X \rangle: \qquad \diag(X) =  \textbf{1}.
\]

Denote by $v_{MC}$ the optimal value of the QUBO, and by $v_{SR}$ the optimal value of the relaxation. Since ${\cal SR}$ is an outer approximation, we obtain the inequality
\[ v_{SR} \leq v_{MC}. 
\]

\medskip

The question now arises how to reconstruct a possibly sub-optimal solution of the QUBO from the optimal solution $X^*$ of the relaxation. Goemans and Williamson proposed the so-called \textit{GW randomized procedure}.

Interpret $X^*$ as the covariance of a Gaussian zero mean vector $\xi \sim {\mathcal N}(0,X^*)$ and define
\[ x = \operatorname{sgn}\xi.
\]
Then $x \in \{-1,1\}^n$ is a random vertex of the hypercube, and defines also a random vertex of the MaxCut polytope ${\cal MC}$.

These generated vertices can be used as sub-optimal solutions of the QUBO. Clearly the expectation of the objective function on these random vertices is an upper bound on the optimal value of the QUBO, since it is an upper bound on the objective value which can be achieved with high probability by generating a sufficient number of samples.

The procedure has the following geometric interpretation. Let $X^* = FF^T$ be a factorization of the rank $k$ relaxation solution, where $F \in \mathbb R^{n \times k}$. Since $\diag(X^*) = \mathbf 1$, the rows $f_i$ of  $F$ are located on the unit sphere $S^{k-1} \subset \mathbb R^k$. 

We may generate the random Gaussian vector $\xi$ according to the formula
\[ \xi = F\psi,
\]
where $\psi \sim {\mathcal N}(0,I)$ is standard normally distributed.
The generated random vertices $x = \operatorname{sgn}\xi$ can then be written element-wise as
\[ x_i = \operatorname{sgn}\langle f_i,\psi \rangle.
\]
This formula admits the following interpretation. The isotropic random hyperplane $H_{\psi} = \{ f \mid \langle f,\psi \rangle = 0 \}$ partitions the rows $f_i$ in two subsets, which correspond to the partition defined by the elements $x_i$.

\medskip

\begin{figure}[h]
    \centering
    \includegraphics[width=0.5\linewidth]{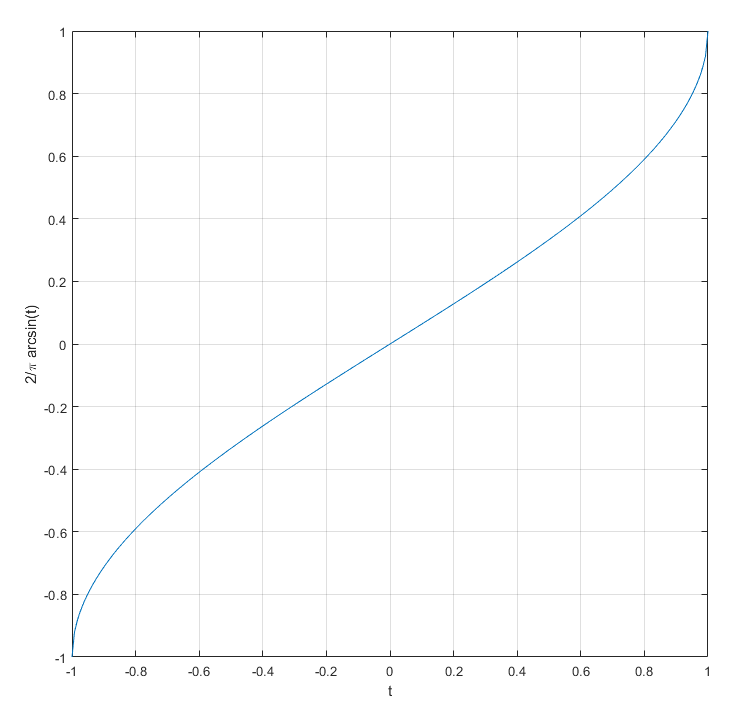}
    \caption{Graph of the function $f(t) = \frac{2}{\pi}\arcsin\,t$ on $[-1,+1]$.}
    \label{fig:arcsin}
\end{figure}

The expectation of the objective is explicitly computable. For $\xi \sim {\cal N}(0,X)$ a Gaussian vector with covariance $X$ satisfying $\diag\,X = {\bf 1}$, it is given by the formula
\begin{equation}\label{expectationFormula}
  \mathbb E_{\xi}(x^TQx) = \frac{2}{\pi}\langle Q,\arcsin X \rangle.
\end{equation} 
Here the function $\arcsin$ is applied element-wise to the matrix $X$. This is possible because the elements of the matrix are contained in the interval $[-1,+1]$. On Fig.~\ref{fig:arcsin} we depict the graph of the function. Note that it is continuously differentiable in the interior of the interval, but the derivative diverges for $X_{ij} \to \pm1$.

\section{Additional Experiments}
\subsection{Effect of Rank on DEM Algorithms}
\begin{figure}[h]
    \centering
    \includegraphics[width=0.9\linewidth]{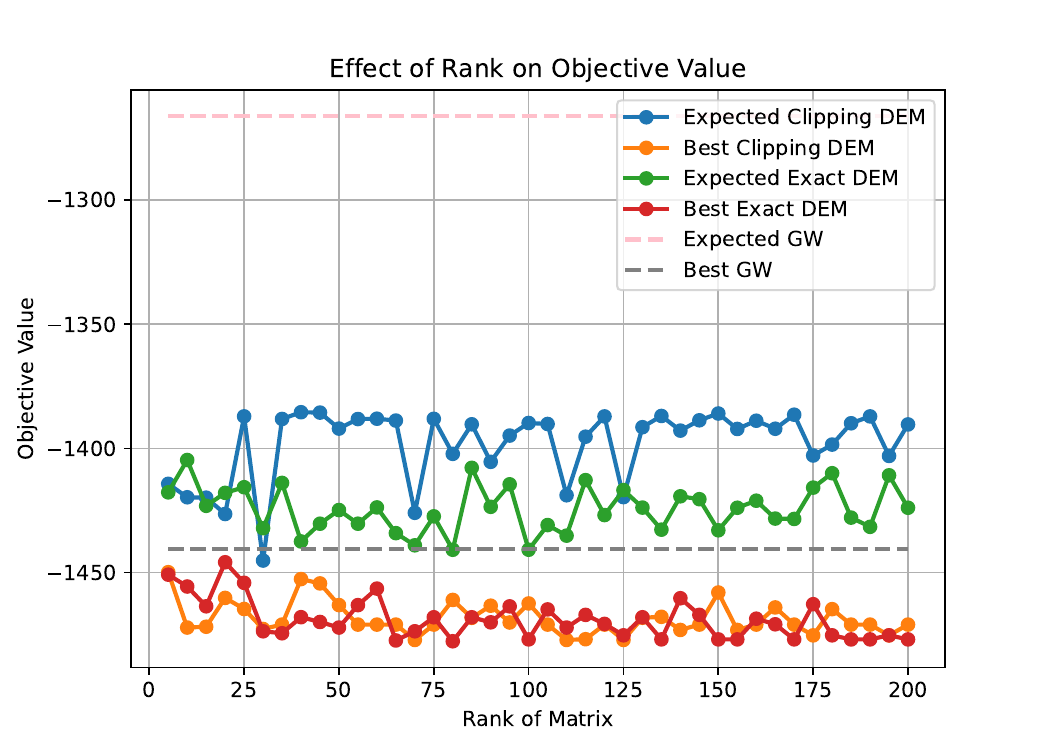}
    \caption{\textbf{Effect of rank on DEM algorithms.} 
    Objective values obtained on a random QUBO instances of size $n=200$ as a function of the rank parameter. 
    Both clipped (DEM-RC) and exact DEM variants are shown, with values for expected and best objectives. 
    For reference, results from the GW-SDP relaxation are included}
    \label{fig:rank_compare}
\end{figure}

We conduct experiments to compare the performance of DEM algorithms across different rank variants. As shown in Figure~\ref{fig:rank_compare}, the Exact DEM method generally achieves lower expected objective values than the clipped variant. However, their best objective values remain nearly identical across most ranks. Notably, the expected objective value of DEM can be reduced to match, and in some cases even surpass, the best value obtained by the GW-SDP relaxation.

\begin{figure}
    \centering
    \includegraphics[width=0.9\linewidth]{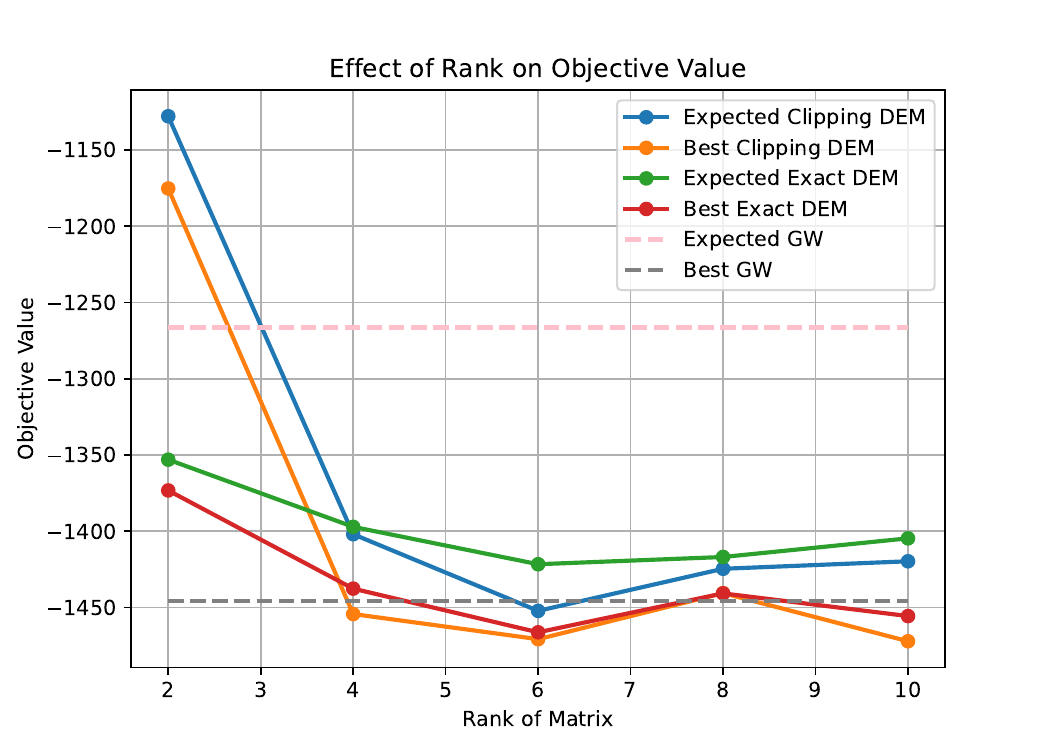}
    \caption{\textbf{Performance at small ranks.}}
    \label{fig:small_ranks}
\end{figure}

We further study the behavior of the DEM at very low ranks. As shown in Figure~\ref{fig:small_ranks}, when the rank is limited to $R=2$, the algorithm performs significantly worse, yielding much higher objective values compared to higher ranks. This suggests that extremely low-rank representations are too limited to describe the structure of the solution space. In contrast, even a small increase in rank (e.g., $R=5,\,10$) is enough to restore competitive performance.

\subsection{Distribution of Objective Values}
\begin{figure}[ht]
    \centering
    \includegraphics[width=0.8\linewidth]{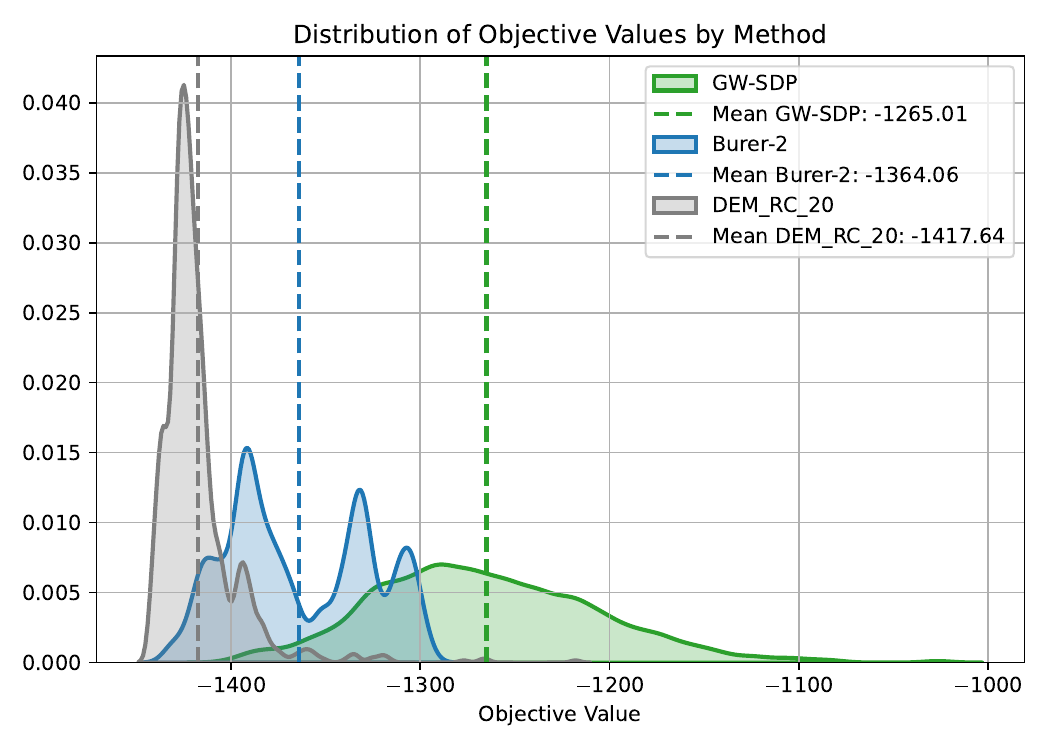}
    \caption{\textbf{Distribution of objective values across methods.}}
    \label{fig:distributions}
\end{figure}

Figure~\ref{fig:distributions} shows the distributions of objective values obtained by the randomized procedure. The DEM method clearly achieves a lower average objective value than the other approaches, highlighting its stronger overall performance. We also notice that methods with worse average results tend to have much higher variability, while DEM produces both better and more consistent solutions across runs.
\subsection{Comparison of the optimization convergence rate of DEM and DEM-RC}
We compare the convergence of DEM and DEM-RC with rank $R=2$. Both methods solve the same optimization problem but use different strategies. DEM applies a DC approach and uses SOCP to compute exact directional derivatives. DEM-RC instead estimates the gradient with a clipping technique to handle non-smooth points at the boundary.

\begin{figure}
    \centering
    \includegraphics[width=0.8\linewidth]{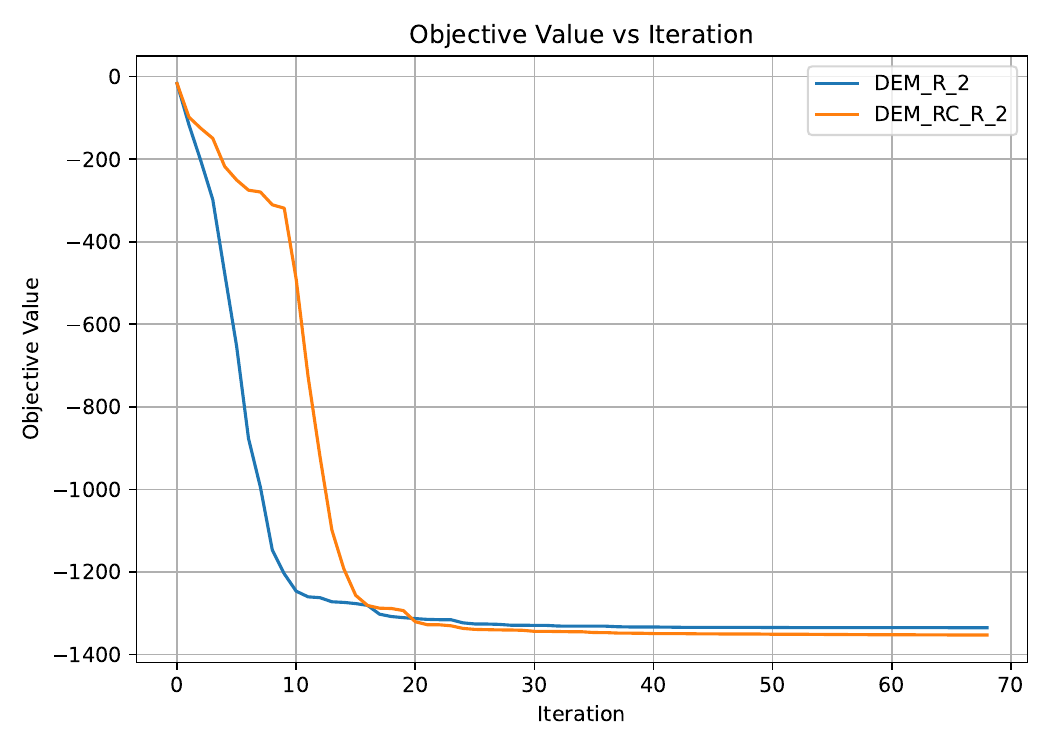}
    \caption{\textbf{Convergence rate of DEM and DEM-RC.}}
    \label{fig:objective_vs_iter}
\end{figure}

Figure~\ref{fig:objective_vs_iter} shows that DEM needs fewer iterations to converge, thanks to its more accurate derivative calculations. However, as seen in Tables~\ref{tab:random_qubo_sizes} and~\ref{tab:IK_results}, DEM-RC runs faster overall, even though it requires more iterations.

\section{Conclusion}

In this work we proposed a method to solve quadratic binary unconstrained optimization problems (QUBO). The solution is generated by the well-known randomized Goemans-Willamson procedure. The covariance matrix of the random Gaussian vector used to generate the vertices of the hypercube is computed in a novel way. 

Instead of relying on the solution of the semi-definite relaxation, we directly minimize the expectation of the objective value on the random vertices. This expectation can be computed analytically, but is not smooth at matrices where an off-diagonal element has absolute value 1. This can be remedied in different ways. Here we use clipping and second-order cone programming. An important detail is that the optimization is not performed directly over the positive semi-definite covariance matrix $X$, which would yield difficult constraints, but over its factor $F$ which is located on a smooth manifold. This gives additional freedom to choose an upper bound on the rank of $X$.

The covariance matrix $\hat X$ obtained by our method yields a better expectation of the objective than the solution $X^*$ of the semi-definite relaxation if we use $X^*$ as a starting point, and it is cheaper to compute. It does not follow that the randomized procedure indeed yields a better solution with $\hat X$ than with $X^*$, but our experiments confirm this hypothesis.

%
%
%
%
\bibliographystyle{splncs04}
\bibliography{icomp2024_conference}

\begin{thebibliography}{10}
\providecommand{\url}[1]{\texttt{#1}}
\providecommand{\urlprefix}{URL }
\providecommand{\doi}[1]{https://doi.org/#1}

\bibitem{ayanzadeh2020reinforcement}
Ayanzadeh, R., Halem, M., Finin, T.: Reinforcement quantum annealing: A hybrid quantum learning automata. Scientific reports  \textbf{10}, ~7952 (2020)

\bibitem{batsheva2023protes}
Batsheva, A., Chertkov, A., Ryzhakov, G., Oseledets, I.: {PROTES}: Probabilistic optimization with tensor sampling. In: NeurIPS 2023 Conference (2023)

\bibitem{Burer02}
Burer, S.: Rank-two relaxation heuristics for max-cut and other binary quadratic programs. SIAM J. Optim.  \textbf{12}(2),  503--521 (2002)

\bibitem{glover1989tabu}
Glover, F.: Tabu search - {P}art {I}. ORSA Journal on Computing  \textbf{1}(3),  190--206 (1989)

\bibitem{GoemansWilliamson}
Goemans, M.X., Williamson, D.P.: Improved approximation algorithms for maximum cut and satisfiability problems using semidefinite programming. J. Assoc. Comput. Mach.  \textbf{42}(6),  1115--1145 (Nov 1995)

\bibitem{goodrich2018optimizing}
Goodrich, T.D., Sullivan, B.D., Humble, T.S.: Optimizing adiabatic quantum program compilation using a graph-theoretic framework. Quantum Information Processing  \textbf{17}(5), ~118 (2018)

\bibitem{HorowitzSahni74}
Horowitz, E., Sahni, S.: Computing partitions with applications to the knapsack problem. J. Assoc. Comp. Mach.  \textbf{21}(2),  277--292 (1974)

\bibitem{kadowaki1998quantum}
Kadowaki, T., Nishimori, H.: Quantum annealing in the transverse {I}sing model. Phys. Rev. E  \textbf{58}, ~5355 (1998)

\bibitem{kirkpatrick1983optimization}
Kirkpatrick, S., Gelatt~Jr., C., Vecchi, M.P.: Optimization by simulated annealing. Science  \textbf{220},  671--680 (1983)

\bibitem{morita2008mathematical}
Morita, S., Nishimori, H.: Mathematical foundation of quantum annealing. J. Math. Phys.  \textbf{49},  125210 (2008)

\bibitem{Nesterov98}
Nesterov, Y.: Semidefinite relaxation and nonconvex quadratic optimization. Optim. Methods Softw.  \textbf{9}(1--3),  141--160 (1998)

\bibitem{salloum2025quantum}
Salloum, H., Savin, S., Kholodov, Y., Ryzhakov, G., Farina, M., Oseledets, I.: Quantum annealing for inverse kinematics in robotics (2025)

\end{thebibliography}


\begin{thebibliography}{8}
\bibitem{ref_article1}
Author, F.: Article title. Journal \textbf{2}(5), 99--110 (2016)

\bibitem{ref_lncs1}
Author, F., Author, S.: Title of a proceedings paper. In: Editor,
F., Editor, S. (eds.) CONFERENCE 2016, LNCS, vol. 9999, pp. 1--13.
Springer, Heidelberg (2016). \doi{10.10007/1234567890}

\bibitem{ref_book1}
Author, F., Author, S., Author, T.: Book title. 2nd edn. Publisher,
Location (1999)

\bibitem{ref_proc1}
Author, A.-B.: Contribution title. In: 9th International Proceedings
on Proceedings, pp. 1--2. Publisher, Location (2010)

\bibitem{ref_url1}
LNCS Homepage, \url{http://www.springer.com/lncs}, last accessed 2023/10/25
\end{thebibliography}

\end{document}